\documentclass[11pt]{amsart}
\usepackage{amsmath,amsthm,epsfig,amssymb}
\title{New parametrization of $A^2+B^2+C^2=3D^2$ and Lagrange's four-square theorem}
\author{Eugen J. Ionascu  }
\curraddr{Department of Mathematics\\ Columbus State University\\4225 University Avenue\\
Columbus, GA 31907} \email{math@ejionascu.edu} \subjclass{}
\date{November $23^{rd}$, 2014}
\keywords{equilateral triangles, system of Diophantine equation,
Eisenstein integers, quadratic forms, Lagrange's four square
theorem, quaternion integers }
\begin{document}
\def\sms{\small\scshape}
\baselineskip18pt
\newtheorem{theorem}{\hspace{\parindent}
T{\scriptsize HEOREM}}[section]
\newtheorem{proposition}[theorem]
{\hspace{\parindent }P{\scriptsize ROPOSITION}}
\newtheorem{corollary}[theorem]
{\hspace{\parindent }C{\scriptsize OROLLARY}}
\newtheorem{lemma}[theorem]
{\hspace{\parindent }L{\scriptsize EMMA}}
\newtheorem{definition}[theorem]
{\hspace{\parindent }D{\scriptsize EFINITION}}
\newtheorem{problem}[theorem]
{\hspace{\parindent }P{\scriptsize ROBLEM}}
\newtheorem{conjecture}[theorem]
{\hspace{\parindent }C{\scriptsize ONJECTURE}}
\newtheorem{example}[theorem]
{\hspace{\parindent }E{\scriptsize XAMPLE}}
\newtheorem{remark}[theorem]
{\hspace{\parindent }R{\scriptsize EMARK}}
\renewcommand{\thetheorem}{\arabic{section}.\arabic{theorem}}
\renewcommand{\theenumi}{(\roman{enumi})}
\renewcommand{\labelenumi}{\theenumi}
\newcommand{\Q}{{\mathbb Q}}
\newcommand{\Z}{{\mathbb Z}}
\newcommand{\N}{{\mathbb N}}
\newcommand{\C}{{\mathbb C}}
\newcommand{\R}{{\mathbb R}}
\newcommand{\F}{{\mathbb F}}
\newcommand{\K}{{\mathbb K}}
\newcommand{\D}{{\mathbb D}}
\newcommand{\comments}[1]{}
\def\phi{\varphi}
\def\ra{\rightarrow}
\def\sd{\bigtriangledown}
\def\ac{\mathaccent94}
\def\wi{\sim}
\def\wt{\widetilde}
\def\bb#1{{\Bbb#1}}
\def\bs{\backslash}
\def\cal{\mathcal}
\def\ca#1{{\cal#1}}
\def\Bbb#1{\bf#1}
\def\blacksquare{{\ \vrule height7pt width7pt depth0pt}}
\def\bsq{\blacksquare}
\def\proof{\hspace{\parindent}{P{\scriptsize ROOF}}}
\def\pofthe{P{\scriptsize ROOF OF}
T{\scriptsize HEOREM}\  }
\def\pofle{\hspace{\parindent}P{\scriptsize ROOF OF}
L{\scriptsize EMMA}\  }
\def\pofcor{\hspace{\parindent}P{\scriptsize ROOF OF}
C{\scriptsize ROLLARY}\  }
\def\pofpro{\hspace{\parindent}P{\scriptsize ROOF OF}
P{\scriptsize ROPOSITION}\  }
\def\n{\noindent}
\def\wh{\widehat}
\def\eproof{$\hfill\bsq$\par}
\def\ds{\displaystyle}
\def\du{\overset{\text {\bf .}}{\cup}}
\def\Du{\overset{\text {\bf .}}{\bigcup}}
\def\b{$\blacklozenge$}

\def\eqtr{{\cal E}{\cal T}(\Z) }
\def\eproofi{\bsq}

\begin{abstract} In this paper we provide a new parametrization for the
diophantine equation $A^2+B^2+C^2=3D^2$ and give a series of corollaries.
We discuss some connections with Lagrange's four-square theorem. As applications, we find new parameterizations
of equilateral triangles and regular tetrahedrons having integer coordinates in three dimensions.
\end{abstract} \maketitle
\section{INTRODUCTION}

In \cite{ejieqtriinz4} we conjectured that the solutions of the diophantine equation

\begin{equation}\label{maindeophantineeq}
A^2+B^2+C^2=3D^2,
\end{equation}

\n can be parameterized in the following way:

{\small \begin{equation}\label{parofmain}
\begin{array}{c} A:=x^2+y^2-z^2-t^2+2(yz+yt+xz-xt),\ \ B:=x^2-y^2+z^2-t^2+2(yz+zt+xt-xy), \\ \\
C:= x^2-y^2-z^2+t^2+2(yt+zt+xy-xz), \ \ D:=x^2+y^2+z^2+t^2 ,\ \text{for some}\ x,y,z\ \text{and}\  t\in\mathbb Z.
\end{array}
\end{equation}
}

In this paper we give a constructive proof of this parametrization and also describe a constructive way of finding all primitive solutions of
(\ref{maindeophantineeq}) for every odd $D$ (an even $D$ simplifies the equation by a factor of $4$) independent of this parametrization.
  Our proofs are only based on the classical arithmetic analysis of Gaussian and Eisenstein integers and some elements of the arithmetic of integer quaternions.
The primitive solutions of (\ref{maindeophantineeq}) (i.e. $0<A\le B\le C$, $\gcd(A,B,C)=1$) for $D\in \{1,3,...,23\}$ are included in the next table:

$$\overset{Table \ 1}{\vspace{0.1in}
\begin{tabular}{|c|c|}
  \hline
  D & [A,B,C] \\
   \hline
  1 & [1,1,1]\\
    3 & [1,1,5]\\
      5 & [1,5,7]\\
       7 & [1,5,11]\\
       9 & [1, 11, 11], [5, 7, 13] \\
 11& [1, 1, 19], [5, 7, 17], [5, 13, 13]\\
  \hline
\end{tabular}\ \ \
\begin{tabular}{|c|c|}
  \hline
  D & [A,B,C]\\
   \hline
13& [5, 11, 19], [7, 13, 17]\\
15& [1, 7, 25], [5, 11, 23], [5, 17, 19]\\
17 & [1, 5, 29], [7, 17, 23], [11, 11, 25], [13, 13, 23]\\
19& [1, 11, 31], [5, 23, 23], [11, 11, 29], [13, 17, 25]\\
21& [1, 19, 31], [11, 19, 29], [13, 23, 25]\\
23& [1, 19, 35], [1, 25, 31], [7, 13, 37], [11, 25, 29]\\
\hline
\end{tabular}}
$$

\n Let us mention that it is known what the number of all primitive solutions of  (\ref{maindeophantineeq}) is (see \cite{ch}, \cite{ejirt}),
and it can be calculated with the formulae

\begin{equation}\label{numberofrepr}
\pi\epsilon(D):=\frac{\Lambda(D)+24\Gamma_2(D)}{48},
\end{equation}
\n where
{\small
\begin{equation}\label{numberofrep2x2plusy2}
\Gamma_2(D)=\begin{cases} 0 \ \text{ if $D$ is divisible by a prime
factor of the form $8s+5$ or $8s+7$}, \ s\ge 0,\\ \\

1 \ \text{ if $D$ is $3$}\\ \\

2^{k}\ \
\begin{cases}\text{ where $k$ is the number of
distinct prime factors of $D$ }   \\ \\
\text {of the form $8s+1$ or  $8s+3$}\  (s> 0),\end{cases}
\end{cases}
\end{equation}}

\begin{equation}\label{2007hirschhorn}
\Lambda(D):=8d\underset{p|D, p\
prime}{\prod}\left(1-\frac{(\frac{-3}{p})}{p}\right),
\end{equation}

\n and $(\frac{-3}{p})$ is the Legendre symbol. Let us recall that, if $p>3$ is a prime, then
{\small
\begin{equation}\label{legendresymbol}
\left(\frac{-3}{p}\right)=\begin{cases}
1\ \  {\rm if }\  p\equiv 1\ {\rm or\ } 7 \ {\rm (mod \ 12)} \\ \\
-1\ \  {\rm if }\  p\equiv 5\ {\rm or\ } 11 \ {\rm (mod \ 12)}
 \ \end{cases},
 \left(\frac{3}{p}\right)=\begin{cases} 1\ \  {\rm if }\  p\equiv 1\ {\rm or\ } 11 \ {\rm (mod \ 12)} \\ \\
-1\ \  {\rm if }\  p\equiv 5\ {\rm or\ } 7 \ {\rm (mod \ 12)}
 \ \end{cases}.
\end{equation}}

\n It is easy to show that every primitive solution of (\ref{maindeophantineeq}) must have $A$, $B$, and $C$ of the form $|6k\pm 1|$ for some $k\in \mathbb Z$.
If $D$ is a prime of the form $4s+1$ then we can find primitive solutions for (\ref{maindeophantineeq}) by writing

$$\begin{array}{c}3D^2=D^2+2D^2=(m^2+n^2)^2+2((m^2-n^2)^2+(2mn)^2)=\\ \\ (m^2+n^2)^2+(m^2-n^2-2mn)^2+(m^2-n^2+2mn)^2,\end{array}$$

\n which gives the only primitive solution for $D=5$ (See Table 1), the primitive solution $(7,13,17)$ for $D=13$ and the primitive solution
$(7,17,23)$ for $D=17$, etc.  If $D$ is a prime of the form $8s+1$ or $8s+3$ ($s \in \mathbb Z$), we have $D=m^2+2n^2$ with $gcd(m,n)=1$ and so

$$3D^2=(m^2-2n^2\pm 4mn)^2+2(m^2-2n^2\mp 2mn)^2,$$

\n which in particular, gives solutions for $D\in \{3,11,17,19\}$ in Table 1. So, the case that remains is $D$ a prime of the form $8s-1$, which in view of Legendre and Gauss's Three Square Theorem, it is exactly the limiting situation when one needs at least four squares to write $D$. If $D$ is a prime of the form $6s+1$ then $D=m^2+3n^2$, so it is a sum of four squares. This gives us in particular the following formulae:

$$3D^2=(m^2-n^2\pm 4nm)^2+(m^2-n^2\mp 4nm)^2+(m^2-5n^2)^2,$$

\n which gives other representations in Table 1 for $D\in \{7,13,19\}$. So, one can reduce the problem of finding primitive representations for (\ref{maindeophantineeq}) to primes
of the form $24s-1$, $s\in \mathbb N$. According to the counting (\ref{numberofrepr}) we have exactly $4s$ primitive solutions in this situation.
One important fact about this counting is that we can say something about the number of partitions of $D$ into sums of four squares (see A001156 and A002635), since in some situations the solutions of (\ref{parofmain}) are uniquely determined (up to a change of signs). Conversely, knowing all the partitions of $D$ as a sum of four squares gives all possible solutions of (\ref{maindeophantineeq}). For example, if $D=23$ then we have only one partition of $23=3^2+3^2+2^2+1^2$ which gives all the representations in Table 1 for $D=23$ by changing the signs in
   (\ref{parofmain}).

(Note: Include here the graph about the number of solutions and the number of partitions)

\section{Main results}

In what follows we will make use of the arithmetic theory of the ring of Eisenstein integers $\mathbb Z[\omega]$.
We refer the reader to \cite{mmeisen} for a
recent presentation of these ideas and a classification of all
the factor rings of $\mathbb Z[\omega]$. It is known that  $\mathbb Z[\omega]$
is an Euclidean domain and as a result every non-unit element has
unique factorization into primes. For an element $z=a+b\omega \in \mathbb Z[\omega]$,
$\omega^2+\omega +1=0$, we use the usual norm
$$N(z)=z\overline{z}=(a+b\omega)(a+b\omega^2)=a^2-ab+b^2.$$
We remind the reader that the units in $\mathbb Z[\omega]$ are
$\{\pm 1,\pm \omega, \pm \omega ^2\}$ and the primes in $\mathbb
Z[\omega]$ are $2$, $1+2\omega$ (the ramified one), every rational
prime $q=6k-1$ and every element of the form $a+b\omega$ or
$b+a\omega$ where $p=N(a+b\omega)=N(b+a\omega)=a^2-ab+b^2$ is a
rational prime of form $6\ell+1$, for some integers $k$, $\ell$ and
$a,b\in \mathbb N$, with $a<b$. Is also known that a natural number $M$ can be represented
as $a^2-ab+b^2$ if and only if the rational prime divisors of the form $2$ or $6k-1$ appear to an even exponent
in the usual prime factorization of $M$. If a prime factor of the form $6k-1$ is present then it divides both $a$ and $b$.

Let us begin with a proof of the Conjecture I we made in \cite{ejieqtriinz4}:

\begin{theorem}\label{conjecturepar} Every primitive integer solution of the Diophantine equation  (\ref{maindeophantineeq}), up to a change of signs, is given by
(\ref{parofmain}).  For every integers $x$,$y$,$z$ and  $t$, the numbers $A$, $B$, $C$
and $D$ defined by (\ref{parofmain}) satisfy  (\ref{maindeophantineeq}).
\end{theorem}

\n \proof.  \ It is cumbersome but just an algebraic exercise to
check that (\ref{parofmain}) gives  $A$, $B$, $C$ and $D$ satisfying
(\ref{maindeophantineeq}). For the non-trivial direction,  let us first start with
an identity that is essential in this constructive proof. We observe
that

$$\begin{array}{c}
[3D-(A+B+C)][3D+(A+B+C)]=9D^2-(A+B+C)^2=\\ \\
3(A^2+B^2+C^2)-(A^2+B^2+C^2+2AB+2AC+2BC)=\\ \\
(A-B)^2+(B-C)^2+(C-A)^2.
\end{array}
$$
We let $U:=3D-(A+B+C)$ and $V:=3D+(A+B+C)$ and observe that since
$C-A=(C-B)-(A-B)$, the above equality can be written as

\begin{equation}\label{ess}
UV=2[(A-B)^2-(A-B)(C-B)+(C-B)^2]=2N(A-B+(C-B)\omega).
\end{equation}

Since $A$, $B$, $C$ and $D$ form a primitive solution of (\ref{maindeophantineeq}) we may assume without loss of generality that all are odd numbers.
Since both of the numbers $U$ and $V$ are even we can write $U=2^{s}\widetilde{U}$ and $V=2^{t}\widetilde{V}$, for some odd integers $\widetilde{U}$ and $\widetilde{V}$ and $s,t\in \mathbb N$. Let us observe that we cannot have both $s$ and $t$ more than or equal to $2$. Indeed, if this the case, then $3D=(U+V)/2$ is divisible by $2$ which contradicts our assumption
on $D$. This implies that $s$ or $t$ is equal to one. Without loss of generality let us say that $s=1$. Then, this forces $t$ even as being the highest power of $2$ dividing
$N(A-B+(C-B)\omega)$. A similar  argument can be used with every prime $q$, $q\equiv -1$ (mod 6), which happens to divide $N(A-B+(C-B)\omega)$. Since $q$ divides $UV$, we know it must divide at least one of the factors. We claim that it cannot divide both of them. Indeed, by way of contradiction if $q$ divides $U$ and $V$, then it must divide $3D=(U+V)/2$ and $A+B+C=(V-U)/2$.
Since $q$ divides $A-B$ and $C-B$ we get it divides $(A+B+C)-(A+C-2B)=3B$. Then, $q$ divides $A$ and $C$ which contradicts our assumption that $gcd(A,B,C)=1$.
Therefore, $q^\alpha$ which appears in the factorization of $N(A-B+(C-B)\omega)$ enters into the decomposition of $U$ or $V$. This shows that $U/2$ and $V$ can be represented in the form
$a^2-ab+b^2$. We will let $\eta=A-B+(C-B)\omega) $.

Let us observe that if the formulae  (\ref{parofmain}) are to be satisfied for some $x$, $y$, $z$ and $t$ we can express $U$ and $V$ in the following way:

$$\begin{array}{c} U=4xt+4yt+4xy+2x^2+2y^2+6z^2+2t^2=2[3z^2+(x+y+t)^2]=\\ \\
=2[(2z)^2-2z(z-x-y-t)+(z-x-y-t)^2]=2N(2z,z-x-y-t),       \ \ \text{ and}
\end{array}$$

$$\begin{array}{c} V=-4xt-4yt-4xy+4x^2+4y^2+4t^2=2[(x-t)^2+(y-t)^2+(x-y)^2]=\\ \\
=4[(x-t)^2-(x-t)(y-t)+(y-t)^2]=4N(x-t,y-t).
\end{array}$$

These calculations suggest that once we have some decompositions for $U/2$ and $V$ in the form $a^2-ab+b^2$ we can pretty much identify
the numbers $2z$, $z-x-y-t$, $x-t$ and $y-t$. This information becomes enough to determine $x$, $y$, $z$ and $t$. However, in order to show that (\ref{parofmain})
are satisfied we need to have ``good" representations for $U/2$ and $V$, which requires a strong relation with $A$, $B$ and $C$ (clearly, we only need the first three equalities in
(\ref{parofmain}) to be satisfied).

Let us denote by $W$ the greatest common divisor of $U/2$ and $V$. We have seen that $W$ must be odd and any prime that divides $W$ must be either $3$ or of the form $6k+1$. Then we can write $U/2=WU'$ and $V=WV'$ with $\gcd(U',V')=1$.

  We define these representations in the following way. We let  $\eta=A-B+(C-B)\omega) $ and first we choose $W'$ and $W''$ such that
 $W'\overline{W'}=W''\overline{W''}=W$ and $W'W''=gcd(W,\eta)$. Of course, this can be done in various ways especially if $W$ has many divisors. Next, we use $W'$ and $W''$ to define
\begin{equation}\label{definitionimportant}
u_1+u_2\omega=W'\gcd(U',\eta/W'W'')\  \ \text{and}\ \ \ v_1+v_2\omega=W''\gcd(V',\eta/W'W'').
\end{equation}

We observe that if $W=1$, these decompositions are unique, up to units in $\mathbb Z[\omega]$. This happens for most of the values of $D$.
Because $$(U/2)V=U'WV'W=U'V'W^2=\eta \overline{\eta} =(\frac{\eta}{W'W''}) W'W'' \overline{\eta}$$
 \n  we can conclude that
\begin{equation}\label{veryimportant} (u_1+u_2\omega)(v_1+v_2\omega)=A-B+(C-B)\omega.
 \end{equation}
\n Because $U/2$ and $V$ have a decomposition as $a^2-ab+b^2$ each of their prime factors must be factors of  $\eta$ or $\overline{\eta}$.
We also notice by taking conjugates, that (\ref{definitionimportant}) implies

\begin{equation}\label{veryimportant1}
U/2=(u_1+u_2\omega)(u_1+u_2\omega^2)=N(u_1,u_2)\ \ \text{ and}\ \ V=(v_1+v_2\omega)(v_1+v_2\omega^2)=N(v_1,v_2).
\end{equation}

We observe that if an odd number can be written as $a^2-ab+b^2$ if both $a$ and $b$ are odd we can use a different representation
in which one of the numbers is even: $a^2-ab+b^2=(b-a)^2-(b-a)b+b^2$. Since $U/2$ is odd we see that this trick allows us to solve the system
$2z=u_1$ and $z-x-y-t=u_2$ in integers. Because $A-B$ and $C-B$ are both even we see that $v_1$ and $v_2$ must be even and so we can
impose that $2(x-t)=v_1$ and $2(y-t)=v_2$. This defines $z=u_1/2$, and then

\begin{equation}\label{problemwith3}
t=\frac{1}{3}[(u_1/2-u_2)-(v_1+v_2)/2]
\end{equation}

\n which apparently is not an integer. Then we solve for $x=v_1/2+t$ and $y=v_2/2+t$.
With these we defined this way we have

$$\begin{cases} A+B+C=(V-U)/2=(v_1^2-v_1v_2+v_2^2)/2-(u_1^2-u_1u_2+u_2^2)\\ \\
A-B=u_1v_1-u_2v_2\\ \\
C-B=u_2v_1+u_1v_2-u_2v_2.
\end{cases}
$$
\n or in terms of $x$, $y$, $z$ and $t$ defined earlier becomes

$$\begin{cases} A+B+C=x^2+y^2-4xy-3z^2-4xt+t^2-4ty\\ \\
A-B=2y^2-2t^2+4zx-2zt-2zy+2xy-2xt\\ \\
C-B=2y^2-2x^2+2zx-4zt+2ty+2zy-2xt.
\end{cases}
$$
Solving this system for $A$, $B$ and $C$ gives exactly (\ref{parofmain}).
So, we only need to show that $t$ is an integer. This is equivalent to showing that $u_1-2u_2\equiv v_1+v_2$ (mod 3).
We know that $V=v_1^2-v_1v_2+v_2^2$ and so $V\equiv 4V\equiv (2v_1-v_2)^2+3v_2^2\equiv(v_1+v_2)^2$ (mod 3).
Similarly $U/2=\equiv (u_1-2u_2)^2$ (mod 3). But $U+V\equiv 0$ (mod 3), which implies $U\equiv 2V$ (mod 3).
We conclude that $u_1-2u_2\equiv \pm (v_1+v_2)$ (mod 3). A change of sign in (\ref{veryimportant}) will take care of this problem.\qed

\vspace{0.1in}
\n {\bf Remark 1:} The following particular case of (\ref{maindeophantineeq}), $A=727$, $B=36293$, $C=85445$ and $D=7(13)(19)31$, has the property
\begin{equation}\label{gcdproperty}
\begin{array}{c} \gcd(D,A+B+C)>1,\  \gcd(D,A+B-C)>1,\\ \\ \gcd(D,A-B+C)>1,\ \ \text{and}\ \  \gcd(D,-A+B+C)>1.\end{array}
\end{equation}
This shows that the case in which $\gcd(U/2,V)>1$ cannot be eliminated in the previous proof.

Next, we exemplify this method. The first case that is non-trivial is $D=3$, $A=1$, $B=1$ and $C=5$.
 Then we get $U/2=(3D-(A+B+C))=1$ and $V=3D+(A+B+C)=16$. Since $gcd(U/2,V)=1$ and $\eta:=A-B+(C-B)\omega=4\omega$, (\ref{definitionimportant})
becomes $u_1+u_2\omega=gcd(1,4\omega)=1$ and $v_1+v_2\omega=gcd(16,4\omega)=4$. Then, we can take  the decompositions
$U/2=(0+\omega)(0+\omega^2)$ and $V=(4+0\omega)(4+0\omega^2)$. This gives $z=0$, $z-x-y-t=1$, $2(x-t)=4$, $2(y-t)=0$. Hence, $x=t+2$, $y=t$ and $x+y+t=-1$. This gives
$t=-1$, $x=1$, $y=-1$ and we obtain back $A$, $B$, $C$ and $D$ with (\ref{parofmain}).

 Let us take an example where $gcd(U/2,V)>1$, $D=49$, $A=5$, $B=17$ and $C=83$.
We get $U/2=21$, $V=252$ and $W=gcd(U/2,V)=21$. Since $$\eta=A-B+(C-B)\omega=6(-2+11\omega)=6(1+2\omega)(1-2\omega)(2+3\omega),$$
\n and $\gcd(W,\eta)=3(2+3\omega)(1-2\omega)$ we can take $$W'=(1+2\omega)(2+3\omega)=-4+\omega\ \text{and}\ \  W''=(1+2\omega)(1-2\omega)=5+4\omega.$$
Then $u_1+u_2\omega=-4+\omega$ and $v_1+v_2\omega=(5+4\omega)\gcd(12,-2(1+2\omega))=2(5+4\omega)(1+2\omega)=6-12\omega$. One can check that
(\ref{veryimportant}) and (\ref{veryimportant1})  are satisfied. Then we obtain $z=-2$, $z-x-y-t=1$, $2(x-t)=6$ and $2(y-t)=-12$.
This gives $t=0$, $x=3$, $y=-6$  and we obtain back $A$, $B$, $C$ and $D$ with (\ref{parofmain}).

\n {\bf Remark 2:}  The parametrization \cite{eji} (page 6) is the particular case of (\ref{parofmain}) for $t=0$ and it was obtained with the chord method.

\subsection{Some consequences and interpretations}

We remind the reader that the Hamilton quaternion algebra over the real numbers, denoted by ${\mathbb H}(\mathbb R)$,
is the associative unital algebra given by requirements:

(i) ${\mathbb H}(\mathbb R)$ is the free $\mathbb R$-module over the symbols $i$, $j$, and $k$, with $1$ the multiplicative unit;

(ii) $i^2=j^2=k^2=-1$, $ij=-ji=k$, $jk=-kj=i$ and $ki=-ik=j$.

By ${\mathbb H}(\mathbb Z)$ we denote the subset of quaternions whose components are all integers. We imbed $\mathbb Z^4$ into ${\mathbb H}(\mathbb Z)$ in the natural way:
$(x,y,z,t)\hookrightarrow x+yi+zj+tk$.  If $q=x+yi+zj+tk\in {\mathbb H}(\mathbb R)$ the conjugate of $q$ is $\overline{q}=x-yi-zj-tk$ and the norm of $q$ is $N(q)=x^2+y^2+z^2+t^2$.
It is known that this norm is multiplicative, i.e. $N(q_1q_2)=N(q_1)N(q_2)$. For some of the elementary results about the arithmetic of ${\mathbb H}(\mathbb Z)$ we recommend the reader the
recent treatment in \cite{dsv}. Given $A$, $B$, $C$ and $D$ satisfying   (\ref{maindeophantineeq}) then the parametrization (\ref{parofmain}) can be written as

\begin{equation}\label{quaternions} Ai+Bj+Ck=q(i+j+k)\overline{q}, \ \text{where}\ \ q=x+yi+zj+tk.
\end{equation}

Let us observe that (\ref{quaternions}) implies  (\ref{maindeophantineeq}) with $D=N(q)$. Our interest in equation (\ref{maindeophantineeq}) came from
finding all of the equilateral triangle having integer coordinates. In the three dimensional space, it turns out that for every solution of  (\ref{maindeophantineeq}) there exists a family of
such equilateral $OP'Q'$ ($O$ the origin) all with the vertices in a plane orthogonal to the vector $(A,B,C)$ determined by a minimal triangle $OPQ$ by formula

\begin{equation}\label{vectorid2}
\overrightarrow{OP'}=m\overrightarrow{\zeta}-n\overrightarrow{\eta},\
\
\overrightarrow{OQ'}=n\overrightarrow{\zeta}+(m-n)\overrightarrow{\eta}
\end{equation}
\n where $\overrightarrow{\zeta}=\overrightarrow{OP}$, $\overrightarrow{\eta}=\overrightarrow{OQ}$ and $m$ and $n$ are arbitrary integers.
To obtain all the minimal triangles one has to implement some rotations of $OPQ$ which is equivalent to taking $(m,n)$ in (\ref{vectorid2}) inn the set
$$ \{(1,0),(1,1),(0,1),(-1,0),(-1,-1),(0,-1)\},$$

\n which we will refer as {\it basic rotations}.
The side-lengths of  $OPQ$ are equal to $D\sqrt{2}$ and those of $OP'Q'$ are given by $D\sqrt{2(m^2-mn+n^2)}$.

In view of the identity (\ref{quaternions}) let us observe that $m-\frac{n}{2}+\frac{n}{2}(i+j+k)$ commutes with $i+j+k$, and so replacing
$q$ with $q'=q[m-\frac{n}{2}+\frac{n}{2}(i+j+k)]$ in (\ref{quaternions}) gives

$$A'i+B'j+C'k=[(m-\frac{n}{2})^2+\frac{3n^2}{4}]q(i+j+k)\overline{q}=(m^2-mn+n^2)(Ai+Bj+Ck),$$

\n which proves that the transformation $q'=q[m-\frac{n}{2}+\frac{n}{2}(i+j+k)]$   leaves the equation (\ref{maindeophantineeq}) essentially the same (just a multiplicative factor).
With all this introductory concept we first have a new way to construct equilateral triangles with integer coordinates in 3-dimensional space.

\begin{corollary}\label{eqtriparametrization}
Every minimal equilateral $OPQ$ triangle in $\mathbb Z^3$ having the origin as one of its vertices, is given (up to a basic rotation) by
\begin{equation}\label{parofeqtri}
\begin{array}{c} P:= [x^2+y^2-z^2-t^2+2(xt-yz), \ -x^2+y^2-z^2+t^2+2(yz+xt),\\  2(yt-zt-xy-xz)],\  \text {and}
 \\ \\
 Q := [x^2+y^2-z^2-t^2+2(yt+xz), \ 2(yz+zt+xt-xy),\\ x^2-y^2-z^2+t^2+2(yt-xz)],
\end{array}
\end{equation}
\noindent  for some $x$, $y$, $z$ and  $t\in\mathbb Z$. Conversely, the points $P$ and $Q$ defined in (\ref{parofeqtri}) together with the origin, form an equilateral triangle
in $\mathbb Z^3$, for  every $x$, $y$, $z$ and  $t\in\mathbb Z$, of sides-lengths equal to $L=(x^2+y^2+z^2+t^2)\sqrt{2}$.
\end{corollary}

This follows from Theorem 2 in \cite{ejiobando} (page 9) and Theorem~\ref{conjecturepar} by taking an arbitrary equilateral triangle, constructing the equation of the plane containing it as in (\ref{maindeophantineeq}), and then using the parametrization one can check that it satisfies the minimality condition of the specified side-lengths. The second part of the statement is  an algebra exercise. An arbitrary equilateral triangle in $\mathbb Z^3$ can be obtained as in Theorem~2 in \cite{ejiobando} and this parametrization doesn't cover all the equilateral triangles since we can have the side-lengths of the form $D\sqrt{2(m^2-mn+n^2)}$, $m$ and $n\in \mathbb Z$, where $D$ is satisfying (\ref{maindeophantineeq}).

\n {\bf Remark:} The formulae (\ref{parofeqtri}) are equivalent to
\begin{equation}\label{parameqtriinqf}
P=q(i-j)\overline{q}, \ \ Q=q(i-k)\overline{q}, \ {\text where}\ \ q=x+yi+zj+tk,
\end{equation}

\n so let us denote by $e_1=i-j$ and $e_2=i-k$.
In \cite{ejirt}, we show that the only triangles that can be extended to regular tetrahedrons in $\mathbb Z^3$ are the ones for which $m^2-mn+n^2=k^2$, and if $k\equiv 0$ (mod 3) then
one can find two solutions for the fourth point.

This suggests that one can cover these triangles with (\ref{parofeqtri}) and so we  can actually go even further
with a new parametrization for regular tetrahedrons.

\begin{corollary}\label{tetrahedronsparametrization}
Every regular tetrahedron $OPQR$  in $\mathbb Z^3$ having the side lengths $\ell \sqrt{2}$ ($\ell \in \mathbb N$) and the origin as one of its vertices is given by $P$, $Q$
as in (\ref{parofeqtri}) and
\begin{equation}\label{paroftetrahedrons}
\begin{array}{c}
\\ \\ R := [2(xt -yt -xz - yz), \ -x^2+y^2-z^2+t^2+2(xy-tz),\\ \\ -x^2+y^2+z^2-t^2-2(zt+xy)], \ \text{and in addition if } \ell\equiv 0\  (mod\ 3)\ \\ \\
R:=[\frac{4}{3}(x^2+y^2-z^2-t^2)+\frac{2}{3}(yz+ty+xz-xt),\\ \\ \frac{1}{3}(x^2-y^2+z^2-t^2) +\frac{8}{3}(yz+xt)+\frac{2}{3}(tz-xy)\\ \\
\frac{1}{3}(x^2-y^2-z^2+t^2)+\frac{8}{3}(ty-zx)+\frac{2}{3}(tz+xy).
\end{array}
\end{equation}
\noindent    for some $x$, $y$, $z$ and  $t\in\mathbb Z$. Conversely, the points $P$, $Q$ and $R$ defined in (\ref{paroftetrahedrons}) together with the origin, form a regular tetrahedron
in $\mathbb Z^3$ for  every $x$, $y$, $z$ and  $t\in\mathbb Z$, with sides-lengths equal to $L=(x^2+y^2+z^2+t^2)\sqrt{2}$.
\end{corollary}

\n \proof. The last statement is purely computational. Let us start with an arbitrary regular tetrahedron $OABC$ (one its vertices is the origin).
The triangle $OAB$ can be written in terms of a minimal one as in (\ref{vectorid2}) with some $m$ and $n$. The minimal triangle can be parameterized as in (\ref{parofeqtri}).
 Since we know that  $m^2-mn+n^2=k^2$ for some $k\in \mathbb N$. It follows that   $k=u^2-uv+v^2$ for some $u,v\in \mathbb Z$ and so
that $m+n\omega=(u+v\omega)^2$. This means that $m=u^2-v^2$ and $n=2uv-v^2$. We observe that
$$\begin{array}{c}
[u-\frac{v}{2}+\frac{v}{2}(i+j+k)](i-j)[u-\frac{v}{2}-\frac{v}{2}(i+j+k)]=\\ \\ (u^2-v^2,2uv-u^2,-2uv+v^2)=(m,n-m,-n)=m(i-j)+n(j-k)=\\ \\
me_1+n(e_2-e_1)=(m-n)e_1+ne_2.
\end{array}$$

\n and

$$\begin{array}{c}
[u-\frac{v}{2}+\frac{v}{2}(i+j+k)](i-k)[u-\frac{v}{2}-\frac{v}{2}(i+j+k)]=\\ \\ (u^2-2uv,2uv-v^2,v^2-u^2)=(m-n,n,-m)=m(i-k)-n(i-j)=\\ \\
=me_2-ne_1=-ne_1+me_2.\end{array}$$
 These formulae are exactly the same as (\ref{vectorid2}) if we take $\overrightarrow{\zeta}=e_1$ and $\overrightarrow{\eta}=e_2$, and make a substitution
$m\to m-n$ and $n\to -n$. This proves that $OAB$ can be obtained with the parametrization (\ref{parofeqtri}) and then we have only two choices for the fourth point.
The rest of the statement follows then from Theorem~2.2 in \cite{ejirt}.\eproof

\n {\bf Remark:} The formulae (\ref{paroftetrahedrons}) are equivalent to
\begin{equation}\label{paramtrhinqf}
P=q(j-k)\overline{q}, \ \ Q=q(\frac{4}{3}i+\frac{1}{3}j+\frac{1}{3}k)\overline{q}, \ {\text where}\ \ q=x+yi+zj+tk.
\end{equation}

\section{Understanding solutions of (\ref{maindeophantineeq})}

Next, let us show how primitive solutions of (\ref{maindeophantineeq}) can be obtained from a given odd $D$. We recall the
following result we obtained in \cite{ejijeff}.
If we consider the sets ${\cal A}:=\{t\in \mathbb Z| t=3x^2-y^2, \gcd(x,y)=1, x,y \in \mathbb Z\}$, ${\cal B}:=\{t\in \mathbb Z| t=x^2+y^2,\gcd(x,y)=1, x,y \in \mathbb Z\}$ and
${\cal C}:=\{t\in \mathbb Z| t=2(x^2-xy+y^2), \gcd(x,y)=1, x,y \in \mathbb Z\}$ then we actually have an interesting relationship between these sets.

\begin{theorem}\label{trinity} ({\bf \cite{ejijeff}}) For the sets defined above, one has the inclusions
\begin{equation}\label{inclusions} {\cal A}\cap {\cal B}\subsetneqq {\cal C},\ \  {\cal B}\cap {\cal C}
\varsubsetneqq {\cal A}, \ \ \text{and}\  \ {\cal C}\cap {\cal A}\varsubsetneqq {\cal B}.
\end{equation}
\end{theorem}

Hence, taking an odd $D$ we write (\ref{maindeophantineeq}) as $A^2+B^2=3D^2-C^2$. So, the question is to find all $C$'s for which
$3D^2-C^2\in {\cal B}$. This makes $3D^2-C^2\in {\cal C}$ also and so $3D^2-C^2$ is in the intersection ${\cal B}\cap {\cal C}\varsubsetneqq {\cal A}$. This implies that such $A^2+B^2=3D^2-C^2$, with $gcd(A,B)=1$,  is made of all the even numbers which have in their decomposition only primes of the form $12k+1$.

Let us observe that if we start with a prime of the form $12k+1$, say $p=13$, we have two writings of $p$,
$p=a^2+b^2=3^2+2^2=u^2-uv+v^2=1^2-1(4)+4^2=(1+4\omega)(1+4\overline{\omega})$. Since $i=\frac{2\omega+1}{\sqrt{3}}$, we can
write $p=(3+2i)(3-2i)=\frac{(3\sqrt{3}+4\omega+2)(3\sqrt{3}-4\omega-2)}{3}\Rightarrow$ $3p=(3\sqrt{3}+1+1+4\omega)(3\sqrt{3}+1+4\overline{\omega})$.
Hence, we get $3p=(3\sqrt{3}+1)^2+(3\sqrt{3}+1)(1+4\omega+1+4\overline{\omega})+(1+4\omega)(1+4\overline{\omega})$ or
$3p=(3\sqrt{3}+1)^2+(3\sqrt{3}+1)(-2)+p$. Therefore, $2p=(3\sqrt{3}+1)(3\sqrt{3}-1)=3(3^2)-1^2$. Since $2p=2(3^2+2^2)=(3+2)^2+(3-2)^2$, we get
the first non-trivial solution of (\ref{maindeophantineeq}): $1^2+1^2+5^2=3(3^2)$. Naturally, the question that arises is whether one can always obtain solutions
of  (\ref{maindeophantineeq}) in a similar way. We know the answer is yes, at the level of existence, from  Theorem~\ref{trinity}, but it is the direct
calculational procedure that we ask if it is always working for every prime of the form $12k+1$. In this respect we have the following proposition whose proof is purely calculatorial.

\begin{proposition}\label{fromprimes} If $p$ is a prime of the form $12k+1$, $p=a^2+b^2=u^2-uv+v^2$ then
$$2p=\frac{(3av+bv-2bu)^2-3(av+bv-2bu)^2}{3v^2-4b^2}.$$
\end{proposition}

The primes of the form $12k-1$ form here an exceptional set and we will be denoting it by $\cal E:=\{11,23,47,...\}$. We know that for every $p\in \cal E$, $3$ is a quadratic residue so let us also denote by $QR_p(3)$ the set of residues $x$ for which $x^2\equiv 3$ (mod $p$). For instance, $QR_{11}(3)=\{5,6\}$,  $QR_{23}(3)=\{7,16\}$,
$QR_{47}(3)=\{12,35\}$, etc. Let us observe that every primitive solution of (\ref{maindeophantineeq}) must have $A\in \{1,5,7,...,D-1\}$.

\begin{theorem}\label{solutionsformain} The equation (\ref{maindeophantineeq}) has primitive solutions for every $A\in \{1,5,7,...,D-1\} $, $\gcd(D,A)=1$, $A\equiv \pm 1$ (mod 6), which is not a solution of any of the equations

\begin{equation}\label{exceptions}
\begin{cases} if \ \ \  p|D,\ \ \ \  A\equiv 0 \ (mod\ p) ,\\
if \ \ \ gcd(D,p)=1,\ \  AD^{-1} \  \ \ \ (mod\ p) \in  QR_p(3).
\end{cases}
\end{equation}
\n for a prime $p=p_i\in {\cal E}$, such that $p_ip_{i+1}\le (3D^2-1)/2$.
 \end{theorem}

\proof. \ \ If $A$ is of the form $6\ell\pm 1$ and $D$ is odd, we can easily check that $3D^2-A^2$ is of the form $2(12k+1)$. Hence, if $A$ does not satisfy (\ref{exceptions}) then $3D^2-A^2$ has in its prime factorization only a factor of $2$ and primes of the form $12k'+1$. Indeed, we remind the reader that $\left(\frac{3}{p}\right)=1$ iff $p=12k\pm 1$. All prime factors of the form $12k-1$ are excluded as divisors of $3D^2-C^2$ because if $p$ is such a divisor, this attracts $3D^2\equiv A^2$ (mod p). From here, we conclude that if  $p|D$ and then $p|A$ but that is excluded, or $(AD^{-1})^2\equiv 3$ (mod p) which is equivalent to $AD^{-1}$ (mod\ p) $\in  QR_p(3)$. The restriction on the size of  $p$ comes from the fact that $3D^2-A^2$ must have at least two factors of the form $12k-1$ if it has one. If the factor has multiplicity two   it does not affect the solvability of the writing $3D^2-A^2$ as a sum of two squares. For every prime of the form $12k+1$ in the factorization of $3D^2-A^2$ we have a representation as a  sum of two coprime squares and so that gives us a writing $3D^2-A^2=B^2+C^2$, with $gcd(B,C)=1$.
\qed

{\bf Remarks:} We observe that we \underline{may} still have solutions for (\ref{maindeophantineeq}) if $A$ satisfies one of the equations (\ref{exceptions}) for some prime $p$, but $p^{\alpha}|| 3D^2-A^2$ ($\alpha$ the greatest power such that $p^{\alpha}| 3D^2-A^2$) with $\alpha$ even.  In general, the number of restrictions (\ref{exceptions}) is smaller than the number of possible values of $A \in \{1,5,7,...,D-1\}$ satisfying $A\equiv \pm 1$ (mod 6). If we do not have any factors of the form $12k+1$ $(k>0$) we end up with a solution in which $B=C$. Since we are taking $A<D$ we still have a non-trivial representation (i.e. not all equal to $D$) but it may not be a primitive representation.

Let us take a few examples to illustrate  what is our approach and for simplicity we will consider only primes for $D$. First, we look at  $D=23$. Since the product of the first two primes in $\cal E$ greater than $11$ is $23(47)=1081>(3(23)^2-1)/2$ the only prime that we have to worry about in Theorem~\ref{solutionsformain} is
$p=11$. Since $23^{-1}$ (mod 11) is $1$ we see that for  all $A\in \{1,11,13,19\}$ (exceptional set is $\{5, 17 \}$) should lead to primitive solutions.
We see that this is indeed the case as shown in Table 1.

For $D=41$, because  $47(59)=2773>2521=(3(41)^2-1)/2$ we need to take into account the first two primes in Theorem~\ref{solutionsformain}. For $p_1=11$, we have to exclude the solutions of
$A\equiv 5D$ (mod 11) or $A\equiv 6D$ (mod 11). This means $A=11s+7$ or $A=11s+4$ for every $s\ge 0$ have to be excluded. We get the restrictions $\{7,29,37\}$. For $p_2=23$,
we have to exclude the solutions of
$A\equiv 7D$ (mod 23) or $A\equiv 16D$ (mod 11). This means $A=23s+11$ or $A=23s+12$ for every $s\ge 0$ have to be excluded. We get the restrictions $\{11,35\}$.
So for $A\in \{1,5,13,17,19,23,25,31,41\}$ we need to have primitive solutions. Indeed, all the primitive solutions for $D=41$ are
$$[1, 1, 71], [5, 23, 67], [5, 47, 53], [13, 43, 55], [17, 23, 65], [19, 31, 61], [25, 47, 47], [31, 41, 49].$$

Let us point out that the solutions of (\ref{maindeophantineeq}) can be organized as a graph in the following way.
First we observe that we have the identity $(2D-a')^2-3D^2=(D-2a')^2-3a'^2$. Two solutions are connected if one can get from one to the other by using this
identity.  We start with the trivial solution for $D=1$ and generate another solution for $D=3$  by using this identity:
 $1^2-3(1^2)=(2-3)^2-3(1^2)=(1-2(3))^2-3(3^2)\Rightarrow$ $1^2+1^2+5^2=3(3^2)$.
Next we proceed as before and replace $1^2-3(3^2)=(2(3)-5)^2-3(3^2)$ with $(3-10)^2-3(5^2)=7^2-3(5^2)$ which implies $1^2+5^2+7^2=3(5^2)$.
At the same time we can consider $5^2-3(3^2)=(2(3)-11)^2-3(3^2)$ and replace it with $(3-2(11))^2-3(11^2)=19^2-3(11^2)$.
Still, $1^2-3(3^2)=(6-7)^2-3(3^2)=(3-2(7))^2-3(7^2)=11^2-3(7^2)$ which gives $1^2+5^2+11^2=3(7^2)$. In Figure~1, we have included a few of the vertices and edges of this infinite graph.

$$\overset{Figure\ 1}{\epsfig{file=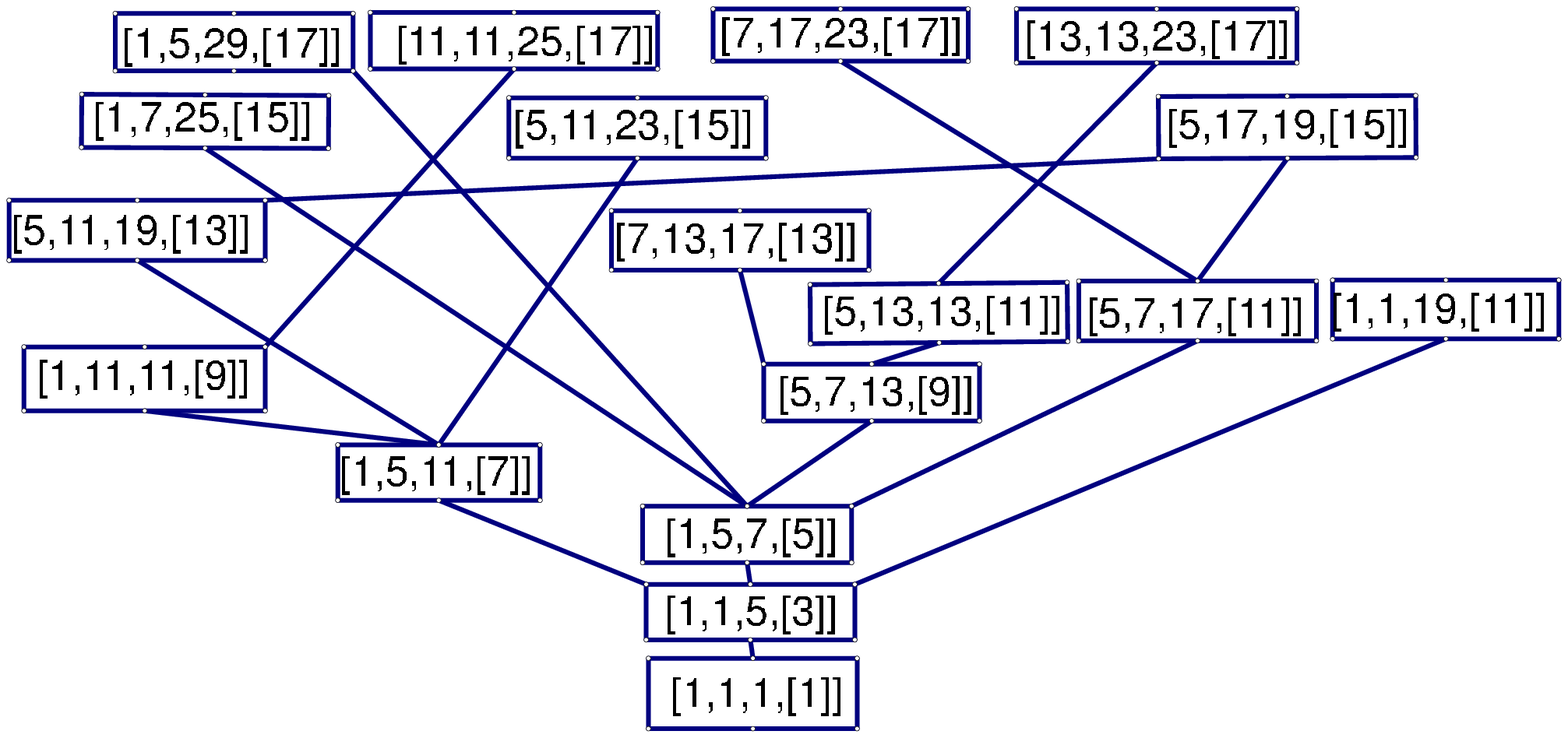,height=2.5in,width=3in} }$$

\n {\bf Acknowledgements:} We thank professor Robin Chapman who shared with us a purely algebraic proof of the parametrization (\ref{parofmain}) using the ($\rm na\ddot{i}ve$) ring theory of integer quaternions.

\end{document}